\documentclass[12pt]{amsart}
\usepackage{geometry}                
\usepackage{graphicx}
\usepackage{amssymb}
\usepackage{epstopdf}
 \usepackage{latexsym}
\usepackage{amsmath}
\usepackage{amsfonts} 
\usepackage{amsthm}

\newtheorem{Proposition}{Proposition}
\newtheorem{Theorem}[Proposition]{Theorem}
\newtheorem{Lemma}[Proposition]{Lemma}

\newtheorem{Corollary}[Proposition]{Corollary}
\newtheorem{Remark}{Remark}

       \def\z{\noindent}

    \def\z{\noindent}  
 
    \def\sqr#1#2{{\vcenter{\vbox{\hrule height .#2pt
                             \hbox{\vrule width .#2pt height#1pt \kern#1pt
                                   \vrule width .#2pt}
                             \hrule height .#2pt}}}}

     \def\CC{\mathbb{C}}
    
    \def\NN{\mathbb{N}}

    \def\ZZ{\mathbb{Z}}

 \def\bfu{\mathbf{u}}
 \def\bfw{\mathbf{w}}
 \def\bff{\mathbf{f}}
  
 \def\bfn{\mathbf{n}}
 \def\bfy{\mathbf{y}}
  \def\bfh{\mathbf{h}}

   \def\bfm{\mathbf{m}}   
      \def\bfF{\mathbf{F}}
\def\bflam{\boldsymbol{\lambda}}
\def\bfphi{\boldsymbol{\phi}}
\def\bfpsi{\boldsymbol{\psi}}

\def\be{\begin{equation}}
\def\ee{\end{equation}}

\begin{document}


\title[Linearization of nonlinear perturbations of Fuchsian systems]{Analytic linearization of nonlinear perturbations of Fuchsian systems}
\author{Rodica D. Costin}




\maketitle

\ 

\ 

\begin{abstract}
Nonlinear perturbation of Fuchsian systems are studied in regions including two singularities. Such systems are not necessarily analytically equivalent to their linear part (they are not linearizable). Nevertheless, it is shown that in the case when the linear part has commuting monodromy, and the eigenvalues have positive real parts, there exists a unique correction function of the nonlinear part so that the corrected system becomes analytically linearizable. 
\end{abstract}

\

\section{Introduction}\label{Introduction}

\subsection{Setting.}\label{Setting}
The paper studies linearization criteria of nonlinear perturbations of Fuchsian systems of the form
\be\label{perFuchs}
\mathcal{E}[\bff]:\ \ \ \ \ \ \ \ \frac{d\bfu}{dx}=M(x)\bfu+\frac{1}{x^2-1}\bff(x,\bfu)\ \ \ \ \ \ \ \ (\bfw\in\CC^d,\ x\in\CC)
\ee
(Systems (\ref{perFuchs}) are denoted by $\mathcal{E}[\bff]$ to distinguish among them by their nonlinear part.)

The linear part 
\be\label{genLinFuchSys}
\mathcal{E}[0]:\ \ \ \ \ \ \ \ \ \ \ \frac{d\bfw}{dx}=M(x)\bfw\ \ \ \ \ \ \ \ (\bfw\in\CC^d,\ x\in\CC)
\ee
is assumed to be a Fuchsian system with three singularities in the extended complex domain.
Their location can be arbitrarily placed using a rational linear transformation, and in the present paper we assume they are located at $x=1,x=-1$ and $\infty$, therefore
\be\label{formM}
M(x)=\frac{1}{x-1}\, A\, +\, \frac{1}{x+1}\, B\ \ \ \ \ \  \ \,A,B\in\mathcal{M}_d(\CC)
\ee

It is assumed that the matrices $A$ and $B$ commute, and formal linearization results are obtained; convergence is proved under the supplementary assumption that the eigenvalues of $A,B$ have positive real parts.

The function $\bff(x,\cdot)$ collects the nonlinear terms in $\bfu$: it has a zero of order two at $\bfu=0$ and 
$\bff(x,\bfu)$ is analytic for $x$ in a simply connected domain $D\ni\{\pm 1\}$ (which will be assumed large enough) and for $\bfu\in\CC^d$, $|\bfu|<R$ (for some $R>0$).

Note that the nonlinear term, $\frac{1}{x^2-1}\bff(x,\bfu)$, is therefore allowed to have at most simple poles at the singularities $x=\pm 1$ of $M(x)$.

 Denote also
\be\label{defQ}
Q(x)=x^2-1
\ee


\subsection{Motivation.}\label{Motivation}

The problem of linearization and the more general, of equivalence, are fundamental in the theory of differential equations, and for many applications and many classes of problems can be reduced to systems of the form (\ref{perFuchs}) -  see \cite{Norm_Form} for references.

The case when the domain $D$ contains just one singular point of $M(x)$ was studied in \cite{Nonln}: generic such systems are linearizable (see \S\ref{rev1s} for details).

When the domain contains two singularities conditions for  formal linearizability of (\ref{perFuchs}), (\ref{formM}) were found in \cite{Norm_Form}, under polynomiality assumptions on the nonlinear part, and not requiring commuting monodromy. 

The present paper shows analytic linearization of canonically corrected systems for general analytic nonlinear part, in the case of commuting monodromy, and positive real parts of eigenvalues of $A$ and $B$.

\subsection{Existence of formal corrections and formal normal forms: results proved in  \cite{Norm_Form}.}\label{Prior_results}

\subsubsection{Assumptions.}\label{eq} \

(a) The function $\bff(x,\bfu)$, which collects the nonlinear terms, has a zero of order two at $\bfu=0$, and 
is analytic for $x$ in the simply connected domain $D\ni\{\pm 1\}$ and for $\bfu\in\CC^d$, $|\bfu|<R$.

(b) The eigenvalues $\boldsymbol{a}=(a_1,\ldots ,a_d)$ of $A$, and respectively $\boldsymbol{b}=(b_1,\ldots ,b_d)$ of $B$, satisfy the Diophantine condition: there exist $C,\nu>0$ so that  
\be\label{DioCond}
\Big| \bfn\cdot\boldsymbol{a} +l-a_s\Big|>C\left( |\bfn|+|l|\right)^{-\nu}\ \ \ {\mbox{and}}\ \ \ \Big| \bfn\cdot\boldsymbol{b} +l-b_s\Big|>C\left( |\bfn|+|l|\right)^{-\nu}
\ee
for all  $l\in\NN$,  $s\in\{1,...,d\}$, and 
$\mathbf{n}\in\NN^d$ with $ |\mathbf{n}|\geq 2$, with the notation
$$\bfn=(n_1\ldots n_d),\ \ \ \ \ \bfn\cdot\boldsymbol{a}=n_1a_1+\ldots+n_da_d,\ \ \ \ \ |\bfn|=n_1+\ldots+n_d$$

(c) No eigenvalue of $A$, or of $B$, is an integer: $a_i,b_i\not\in\ZZ,\ i=1,\ldots ,d$.

(d) The eigenvalues $\lambda_1,\ldots ,\lambda_d$ of the matrix $A+B$ satisfy the following nonresonance condition:
\be\label{non_res}
k+\mathbf{n}\cdot {\bflam}-\lambda_j \not=0\ \ \ {\mbox{for\ all}}\ \ \mathbf{n}\in\NN^d\ ,\ k\in\NN\ ,\ \ j=1,\ldots ,d
\ee

\subsubsection{Obstructions to linearization and formal correction.}\label{formal}

The following result was proved in  \cite{Norm_Form}. Note that there are no commutativity assumptions on $A$ and $B$.

\begin{Theorem}\label{Obstructions}

Consider the system (\ref{perFuchs}) under the assumptions of \S\ref{eq}.

Assume that the nonlinear part is of polynomial type in $x$, in the sense that
\be\label{pol_type}
{\mathbf{f}}(x,\bfu)=\sum_{|\mathbf{m}|\geq 2}\mathbf{f}_\mathbf{m}(x)\mathbf{u}^\mathbf{m}\ \ \ \ \ {\mbox{with}}\ \mathbf{f}_\mathbf{m}(x)\ {\mbox{polynomials}}
\ee
with the usual notation
$$\bfu^\bfm=u_1^{m_1}u_2^{m_2}\ldots u_d^{m_d}\ \ \ {\mbox{for\ }} \bfu=(u_1,u_2,\ldots ,u_d)\in\CC^d,\ \bfm\in\NN^d$$

Then {\em{there exists a unique}} correction $\bfphi(\bfu)$ of $\bff(x,\bfu)$ as a formal series 
\be\label{serphi}
\bfphi(\bfu)=\sum_{|\bfn|\geq 2} \bfphi_\bfn \bfu^\bfn
\ee
so that the correction  $\mathcal{E}[\bff-\bfphi]$ of $\mathcal{E}[\bff]$:
\be\label{corperFuchs}
\mathcal{E}[\bff-\bfphi]:\ \ \ \ \ \frac{d\bfu}{dx}=\left(\frac{1}{x-1}\, A\, +\, \frac{1}{x+1}\, B\right)\bfu+\frac{1}{x^2-1}\left[\bff(x,\bfu)-\bfphi(\bfu)\right]
\ee
is linearizable by a formal series transformation
\be\label{serh}
{\bfu=\mathbf{H}(x,\bfw)=\bfw+\sum_{|\mathbf{m}|\geq 2}\mathbf{h}_\bfm(x)\mathbf{w}^\mathbf{m}}
\ee
where $\mathbf{h}_\bfm(x)$ are functions analytic on $D$ (in fact, they are polynomials).

\end{Theorem}

Note that the coefficients $\bfphi_\bfn$, $\bfn\in\NN^d, |\bfn|\geq 2$ represent the obstructions to linearization, in the sense that the equation $\mathcal{E}[\bff]$ is linearizable if and only if all $\bfphi_\bfn$ are zero.

\subsubsection{Remark.} Obviously, formal linearization is a necessary condition for analytic linearization.

\section{Main Results.}\label{Main_res}

While the assumptions of \S\ref{eq} are essential for the results of Theorem\,\ref{Obstructions} to hold, it is natural to expect that  the condition (\ref{pol_type}) that $\bff$ be of polynomial type in $x$ to be generalizable to holomorphic functions of $x$. It is also to be expected that the series (\ref{serphi}) and (\ref{serh})  converge if the eigenvalues of $A+B$ are not "too close to resonance" (see (\ref{non_res})).

The present paper considers the case when the matrices $A$ and $B$ commute, and when their eigenvalues have positive real parts. Under these assumptions there are no polynomiality assumptions on the nonlinear part $\bff$, which can be any holomorphic function (on a domain large enough). It is shown that a unique formal correction $\bfphi(\bfu)$ still exist {{for functions $\bff(x,\bfu)$ analytic in $x$}}  - Theorem \ref{anf}, and furthermore,  that the linearization series (\ref{serh}), and the series 
(\ref{serphi}), converge, hence we have {\em{analytic linearization }}- Theorem \ref{ConvSer}.

\subsection{Setting.}\label{eq1}

Consider the equation $\mathcal{E}[\bff]$ given by (\ref{perFuchs}),(\ref{formM}) 
under the assumptions of \S\ref{eq}.

Furthermore, it is assumed that $A$ and $B$ are simultaneously diagonalizable:
\be\label{ass_diag}
A={\mbox{diag}} (a_1,\ldots ,a_d),\ \ \ \ \ B={\mbox{diag}} (b_1,\ldots ,b_d)
\ee

Theorem\,\ref{ConvSer} is proved under the additional assumption
\be\label{rhp}
\Re a_j>0,\ \ \ \ \ \Re b_j>0\ \ \ \ \ {\mbox{for\ all\ }}j=1,\ldots, d
\ee

\subsection{Obstructions to linearization - equations with analytic nonlinear terms.}\label{SecObs}

\begin{Theorem}\label{anf}

Consider equation $\mathcal{E}[\bff]$  given by (\ref{perFuchs}),(\ref{formM}) under the assumptions of \S\ref{eq1}.
Let $\bff$ be holomorphic for $x\in D$ and $|\bfu|<R$.

Then{\em{ there exists a unique}} "correction" $\bfphi(\bfu)$ of $\bff(x,\bfu)$ as a formal series (\ref{serphi})
so that the "corrected" equation $\mathcal{E}[\bff-\bfphi]$ - see (\ref{corperFuchs}) - 
is linearizable by a formal series (\ref{serh}) with coefficients $\bfh_\bfm(x)$ holomorphic on $D$.

\end{Theorem}

\subsection{Analytic linearization.}

Theorem\,\ref{ConvSer} shows that the series in Theorem\,\ref{anf} converge, therefore any nonlinear term $\bff$ can be corrected to an analytic function so that the resulting equation is analytically linearizable.

\subsubsection{{The numbers $c$ and $\rho_{min}$.}}\label{rho_min}

These quantities are used in the statement of Theorem\,\ref{ConvSer} to define a domain $D$ (where the nonlinear term will be required to be analytic) and are defined as follows.  Denote
\be\label{cm}
c_\bfm=\frac{(\bfm/n)\cdot(\mathbf{b}-\mathbf{a})}{(\bfm/n)\cdot(\mathbf{b}+\mathbf{a})}\ \ \ \ \ \ \  {\mbox{for}}\ \ \ \bfm\in\NN^d,\ n=|\bfm|\geq 2
\ee
 
The set of all points $c_\bfm$ given by (\ref{cm}) belong to a compact set $K$:
$$K=\left\{\, g(\mathbf{t})=\frac{\mathbf{t}\cdot(\mathbf{b}-\mathbf{a})}{\mathbf{t}\cdot(\mathbf{b}+\mathbf{a})}\, ;\, \mathbf{t}\in [0,1]^d,\, |\mathbf{t}|=1\,\right\}$$
(In fact the boundary of $K$ consists of arcs of the circles $g(\mathbf{t})$ obtained for $\mathbf{t}$ having only two nonzero coordinates.)

Let $c\in K$ and let $\rho_{min}$ positive and large enough so that the disk $|x-c|<\rho_{min}$ contains the points $-1$, $1$ and the set $K$.

\subsubsection{Analytic linearization theorem}
\begin{Theorem}\label{ConvSer}

Consider the system (\ref{perFuchs}),(\ref{formM}) under the assumptions of \S\ref{eq1}.
Let $c$ and $\rho_{min}$ as in \S\ref{rho_min}.

Assume that the nonlinear part $\bff(x,\bfu)$ is holomorphic for $|\bfu|<R_0$ and $x$ in a disk
$|x-c|<\rho_0$ where $\rho_0>\rho_{min}$.

Then the series (\ref{serphi}) and (\ref{serh}) converge in a subdomain $|\bfu|<R_e$ and $|x-c|<\rho_e$ ($\rho_{min}<\rho_e<\rho_0$, $0<R_e<R_0$). 

Furthermore, $\rho_e$ can be made arbitrarily close to $\rho_0$ if $R_e$ is small enough.

\end{Theorem}

\subsubsection{Comments regarding the proofs.}

{\em{(i)}} The proof of Theorem\,\ref{Obstructions} in \cite{Norm_Form}, which is done under polynomiality assumptions of the nonlinear terms, relies on  an algebraic structure: the terms are expanded  in a base of polynomials in $x$ satisfying a generalized Rodrigues formula. Under the additional commutativity assumption of the present paper, these special polynomials are Jacobi polynomials (see also \cite{Gen_Jacobi}). 

The natural approach to proving {Theorem}\,\ref{ConvSer} seems to be by expansions in Jacobi series in the variable 
$x$, on an elliptic domain $D$ with foci at $x=\pm 1$, see \cite{Carlson}, and the condition (\ref{rhp}) could be replaced by a Diophantine condition in (\ref{non_res}).
However, some technical results needed for estimates are not available yet, or at least, are not known to the author. The proof of convergence is done here using estimates on a disk $D$ by saddle-point and Laplace's method, \S\ref{PhiJ}. To avoid small denominator problems a rapidly convergent iteration is used in the proof. The estimates of \S\ref{homeq}-\S\ref{convseqem} follow the line in \cite{Arnold} \S 28. Additional refinements were needed in the proof of convergence of the correction (\ref{serphi}) in \S\ref{convcorr}.

{\em{(ii)}} The fact that a formally corrected system which is formally linearizable turns out to be analytically linearizable under the positivity assumption (\ref{rhp}) is a result in the spirit of the Poincar\'e-Dulac theorem (see e.g.  \cite{Arnold}). It should be noted that this theorem cannot be used in the present context because after subtracting the formal correction $\boldsymbol{\phi}(\bfu)$ the equation is only presented as a formal series, hence not known to be analytic. The convergence of both the linearization map and of the correction series is proved here simultaneously.


\section{Proofs.}\label{Proofs}

\subsection{Proof of Theorem \ref{anf}.}\label{PfOb}

\subsubsection{Obstructions to linearization.}\label{recsys}

Section \S\ref{recsys} contains the first steps in the proof. They follow \cite{Norm_Form} and are reproduced here for completeness.

If a map (\ref{serh}), $\bfu=\mathbf{H}(x,\bfw)=\bfw+\bfh(x,\bfw)$, is a linearization map of (\ref{perFuchs}) then then it satisfies the nonlinear partial differential equation
\be\label{nonlinPDE}
\partial_x\bfh+d_\bfw\bfh\, M\bfw-M\bfh=\frac{1}{x^2-1}\, \left[\bff(x,\bfw+\bfh)-\bfphi(\bfw+\bfh)\right]
\ee
where $M$ is given by (\ref{formM}).

Denote by $\bfh_n$ the homogeneous part of degree $n$ of the function $\bfh(x,\bfw)$
\be
\bfh_n(x,\bfw)=\sum_{|\mathbf{m}|=n}\bfh_\mathbf{m}(x)\bfw^\mathbf{m}
\ee
and use a similar notation for other functions ($\bff(x,\bfw)$ etc.).
Expanding in power series in $\bfw$ equation (\ref{nonlinPDE}) yields a recursive system for $\bfh_n$:
\be\label{eqhn}
\partial_x\bfh_n+d_\bfw\bfh_n\, M\bfw-M\bfh_n=\frac{1}{x^2-1}\, \mathbf{R}_n(x,\bfw)\ \ \ \ \ (n\geq 2)
\ee
where
\be\label{formRn}
\mathbf{R}_n=\bff_n-\bfphi_n+\tilde{\mathbf{R}}_n
\ee
with $\tilde{\mathbf{R}}_n$ a
polynomial in $\bfphi_\mathbf{m}$, $\bfh_\bfm$, $\bff_\mathbf{m}$ with $|\mathbf{m}|<n$.

Denote by $Y(x)$ a fundamental matrix for the linear part: $Y'=MY$. Using the variation of parameters formula for a linear nonhomogeneous equation, and choosing the solution which is not branched at $x=-1$ we obtain
\be\label{hnat-1}
\bfh_n(x,\bfw)=Y(x)\int_{-1}^xQ(t)^{-1}Y(t)^{-1}\mathbf{R}_n(t,Y(t)Y(x)^{-1}\bfw)dt
 \ee
is a  particular solution of (\ref{eqhn}), and this solution is analytic at $x=-1$ (see the Appendix, \S\ref{details} for details).

Rewriting (\ref{hnat-1}) as
\be\label{splithn}
\bfh_n(x,\bfw)=Y(x)\int_{-1}^1Q^{-1}Y^{-1}\mathbf{R}_ndt+Y(x)\int_{1}^xQ^{-1}Y^{-1}\mathbf{R}_ndt
 \ee
the last term of (\ref{splithn}) is the  unique linearization map of  (\ref{perFuchs}) which is analytic at $x=1$.

 Then the linearization map $\bfh$ is analytic at both $x=1$ and $x=-1$ if and only if the first term in 
right hand side of (\ref{splithn}) vanishes:
\be\label{obcond}
\int_{-1}^1Q(t)^{-1}Y(t)^{-1}\mathbf{R}_n(t,Y(t)Y(x)^{-1}\bfw)\,  dt\, =\, 0\ \ \ {\mbox{for\ all\ }}\bfw\in\CC^d,\ n\geq 2
\ee

Formulas (\ref{obcond}) show the obstructions to linearization: there is, for every $\mathbf{m}\in\NN^d, |\mathbf{m}|\geq 2$, one numerical condition (vector-valued in $\CC^d$).

\subsubsection{Existence of the correction $\boldsymbol{\phi}(\bfu)$.}\label{extphi}

Using  (\ref{formRn}) in (\ref{obcond}) we obtain recursively for $n\geq 2$ equations for $\bfphi_n(\bfw)$ of the form
\be\label{sysphi1}
\int_{-1}^1Q(t)^{-1}Y(t)^{-1}\bfphi_n(Y(t)Y(x)^{-1}\bfw)\,  dt\, =\bfF_n(x,\bfw)
\ee
where $\bfF_n(x,\bfw)$ are homogeneous polynomials in $\bfw$ degree $n$, with coefficients vector-valued functions, analytic in $x$ on $D$.

The matrices $A$ and $B$ being given by (\ref{ass_diag}), then $Y(x)$ is the diagonal matrix:
$$Y(x)={\mbox{diag}}\, \left[ y_1(x), \ldots , y_d(x)\right]\ \ \ \ {\mbox{where}}\ y_j(x)=(x-1)^{a_j}(x+1)^{b_j}$$
and (\ref{sysphi1}) becomes

\be\label{matsys}
\sum_{|\bfm|=n}\, \frac{1}{\bfy (x)^\bfm}\, \int_{-1}^1\, \, (t-1)^{\bfm\cdot\mathbf{a} -1}(t+1)^{\bfm\cdot\mathbf{b} -1}Y(t)^{-1}\, dt\, \bfphi_\bfm\bfw^\bfm\, =\, \sum_{|\bfm|=n}\, \bfF_\bfm(x)\bfw^\bfm
\ee
Equation (\ref{matsys}) is a linear system for $\{\bfphi_\bfm\}_{|\bfm|=n}$. 
The left-hand-side of (\ref{matsys}) is a diagonal linear operator having as entries Eulerian integrals of the first kind 
$$\frac{1}{\bfy (x)^\bfm}\, \int_{-1}^1\, \, (t-1)^{\alpha_{\bfm,j} -1}(t+1)^{\beta_{\bfm,j} -1}\, dt\, =\, \frac{(-1)^{\alpha_{\bfm,j} -1}2^{\alpha_{\bfm,j}+\beta_{\bfm,j}-1}B(\alpha_{\bfm,j},\beta_{\bfm,j})}{\bfy (x)^\bfm}$$
where $\alpha_{\bfm,j}= \bfm\cdot\mathbf{a}-a_j$,  $\beta_{\bfm,j}= \bfm\cdot\mathbf{b}-b_j$, and  $B(p,q)$ is the beta function
\be\label{Bfun}
B(p,q)=\int_{0}^1s^{p-1}(1-s)^{q-1}ds=\frac{\Gamma(p)\Gamma(q)}{\Gamma(p+q)}
\ee
Therefore equation (\ref{matsys}) has a unique solution if the nonresonance condition (\ref{non_res}) is satisfied, and Theorem \ref{anf} is proved. \qed

\subsection{Proof of {Theorem} \ref{ConvSer}.}

The proof is written in dimension one, for simplicity of notation. The only (minor) difference for dimension $d>1$ is mentioned in  \S\ref{ddim}.


\subsubsection{Initial steps.}

Note that in view of Theorem\,\ref{anf} it can be assumed that $f_n=0$ for $2\leq n< N$ for any chosen $N>2$. It follows that also $\phi_n=0$ and $h_n=0$ for $n<N$.

\subsubsection{The functionals $\Phi_n$ and the operators $J_n$.}\label{PhiJ}

We will consider $\rho$ in a closed interval 
$I=[{\rho}_{min},\rho_{0}] $
and  $\delta\in (0, 1/2)$.

Let $c=\frac{b-a}{b+a}$. Denote $D_\rho=\{x\in\CC\, ;\, |x-c|<\rho\}$. Let $\mathcal{B}_\rho$ be the Banach space of functions analytic on $D_\rho$, continuous on $\overline{D_\rho}$, with the sup norm:
$$\|F\|=\sup_{x\in D_\rho}|F(x)| $$
Denote for simplicity the sup norm on $\mathcal{B}_{\rho\rm{e}^{-\delta}}$ by $ \|\cdot \|_\delta$:
$$ \|F\|_\delta=\sup_{x\in D_{\rho\rm{e}^{-\delta}}}|F(x)|$$

\z Note that using a Cauchy integral formula (see also Remark\,\ref{R6}) we have
\be\label{estimF'}
\|F'\|_\delta \leq \frac{\|F\|}{\rho\left( 1-\rm{e}^{-\delta}\right)}\leq \frac{2}{\rho\delta}\, \|F\|
\ee

Let $n\geq 1$. For  $F\in\mathcal{B}_\rho$ define the linear functionals
\be\label{defPhi}
\Phi_{n+1}[F]=\frac{2^{1-n(a+b)}}{B(na,nb)}\, \int_{-1}^1 (1-t)^{na-1}(1+t)^{nb-1}\, F(t)\, dt
\ee
\z (see (\ref{Bfun})) and the linear operators

\be\label{opJ}
J_{n+1}[F](z)=(1-z)^{-na}(1+z)^{-nb}\, \int_{-1}^z\, (1-t)^{na-1}(1+t)^{nb-1}\, \left( F(t)-\Phi_{n+1}[F]\right)\, dt
\ee
where $\Phi_{n+1}[F]$ is given by (\ref{defPhi}) and therefore $J_{n+1}[F]$ is a function analytic at both endpoints $z=\pm 1$ (see \S\ref{recsys} and \S\ref{extphi}).

Note that
\be\label{invalPJ}
\Phi_n[1]=1\ ,\ J_n[1]=0\ ,\ \Phi_n[x-c]=0\ ,\ J_{n+1}[x-c]=-\frac{1}{n(a+b)}
\ee

Denote 
\be\label{tildenot}
F(x)=\sum_{k\geq 0}F_k(x-c)^k=F_0+(x-c)F_1+\tilde{F}(x)
\ee

\begin{Lemma}\label{LemmaPhi}

If $F\in\mathcal{B}_\rho$ then
$$\Phi_{n+1}[F]=F(c)+\frac{1}{n+1}\tilde{\Phi}_{n+1}[\tilde{F}]$$
with 
\be\label{estimRn}
| \tilde{\Phi}_{n+1}[\tilde{F}] | \leq \, {\mbox{const}}\ \max \left\{\| F\|, \|F'\|\right\}
\ee
where the constant is independent of $\rho\in I$.


\end{Lemma}

\begin{Lemma}\label{LemmaJ}

If $F\in\mathcal{B}_\rho$ then $J_{n+1}[F]\in\mathcal{B}_{\rho{\rm{e}}^{-\delta}}$ and 
$$J_{n+1}[F]=-\frac{1}{n(a+b)}\frac{F(x)-F(c)}{x-c}+\frac{1}{n^2}\tilde{J}_{n+1}[\tilde{F}]$$
with 
\be\label{estimSn}
\Big\| \tilde{J}_n[\tilde{F}](x) \Big\|_{\delta/2}\, \leq \, {\mbox{const}}\ \  \max \left\{ \|F\|,\|F'\|\right\}
\ee
where the constant is independent of $\rho\in I$ and $\delta\in (0,\frac{1}{2})$.

In particular, 
$$\| \tilde{J}_n[\tilde{F}] \|_\delta \, \leq \, {\mbox{const}}\ \ \delta^{-1}\| F\|$$

\end{Lemma}

\begin{Corollary}\label{Estim1}
The operators $J_n$ satisfy the estimates
\be\label{Estim1}
\big\|J_{n+1}[F]\big\|_\delta\leq \ \frac{ {\mbox{const}}}{n}\ \delta^{-1}\|F\|
\ee
with the constant not depending on $\delta$ or $\rho$ for all $\rho\in I$ and $\delta\in (0,\frac{1}{2})$. 
\end{Corollary}

\z {\bf{Proofs}}

The proofs of Lemma\,\ref{LemmaPhi} and Lemma\,\ref{LemmaJ}, respectively,  are done by finding the asymptotic behavior of the integrals for large $n$.

\z{\em{I. The lines of steepest ascent.}}

Consider the function 
\be\label{fung}
g(x)=a\ln(1-x)+b\ln(1+x),\ \ \ \ \ x\in D_\rho
\ee

The function $\Re g(x)$ has a saddle-point at $x=c$. The stable manifold intersects the real line at $c_0=\frac{\Re (b-a)}{\Re(b+a)}$. To the left of the stable manifold steepest descent lines wind toward $x=-1$, while to the right of this manifold they wind towards $x=1$.

A straightforward calculation shows that, for any $z$, $\Re g(x)$ is increasing on the segment $[-1,z]$  (steep ascent) if this segment does not intersect the disk $\mathcal{C}_-$ bounded by the circle passing through the points $1$, $c_0$ and $c$.

Similarly, $\Re g(x)$ increasing on the segment $[1,z]$ (steep ascent) if this segment does not intersect the disk $\mathcal{C}_+$ bounded by the circle through the points $-1$, $c_0$ and $c$.

\z{\em{II. Proof of Lemma \ref{LemmaPhi}.}}
To obtain the behavior of the integral for large $n$ the path of integration in (\ref{defPhi}) must be taken to pass through $c$, along the unstable manifold. The saddle-point method (see, e.g., \cite{Bender-Orszag}) gives
$$\Phi_{n+1}[F]\sim \frac{2^{1-n(a+b)}}{B(na,nb)}\, {\rm{e}}^{ng(c)}\, F(c)\, \int_{-\infty}^{\infty} {\rm{e}}^{ng''(c)/2 \, s}\, ds\ \ \ \ {\mbox{as}} \ n\to\infty$$
and Lemma \ref{LemmaPhi} follows by direct calculation and standard estimates of the remainder.

\z{\em{III. Proof of Lemma \ref{LemmaJ}.}}


It was shown that the maximum modulus of $(1-x)^a(1+x)^b$ is attained for $x=z$ for any $z\not\in \mathcal{C}_-\cap\mathcal{C}_+$. (Note that if $c_0=c$ these disks are tangent. This happens for, e.g.  $a=b$, or $a,b>0$.)

For  $z\in \mathcal{C}_-\cap\mathcal{C}_+$ consider the following path of integration in (\ref{opJ}), described here in opposite direction: the line of steepest descent through $z$ until this line exits $\mathcal{C}_-\cap\mathcal{C}_+$,  followed by a segment to $-1$ (steep descent) if $z$ is on the left-side of the unstable manifold, or on a segment to $1$, if $z$ lies to the right of this manifold. 
Then, also in this case the maximum modulus of $(1-x)^a(1+x)^b$ is attained at $x=z$.

Therefore in all cases the Laplace method yields
$$J_{n+1}[F](z)\, \sim\,  {\rm{e}}^{-ng(z)}\, \int_\cdot^z\, {\rm{e}}^{ng(z)+ng'(z)(t-z)}\, \frac{F(z)-F(c)}{1-z^2}\, dt\ \ \ {\mbox{as\ }}n\to\infty$$
and Lemma \ref{LemmaJ} follows by direct calculation and standard estimates of the remainder. \qed


\subsubsection{The homological equation.}\label{homeq}

The dominant part of (\ref{nonlinPDE}) is

\be\label{homologeq}
\partial_x h+d_w h\, Mw=Mh+\frac{1}{1-x^2}\, \left( f(x,w)-\Phi[f](w)\right)
\ee
where $\Phi[f](w)$ is the unique power series for which equation (\ref{homologeq}) has a solution $h$ as a power series in $w$ with coefficients analytic in $x$:
$$\Phi[f](w)=\sum_{n\geq 2}\, \Phi[f]_n\, w^n\ \ \ {\mbox{where}}\ \ \ \Phi[f]_n=\Phi_n[f_n]$$

If $h$ satisfies (\ref{homologeq})  then the map $u=H(x,w)=w+h(x,w)$ is an equivalence map between the equation
$$u'=M(x)u+\frac{f(x,u)-\Phi[f](u)}{1-x^2}$$
and 
$$w'=M(x)u+\frac{\mathcal{R}[f](x,w)}{1-x^2}$$
where
\be\label{defR}
\mathcal{R}[f](x,w)=\left( {1+d_wh} \right)^{-1}\, {{\mathcal{R}_0}[f](x,w)} 
\ee
with
\be\label{defD}
\mathcal{R}_0[f](x,w)=f(x,w+h)-f(x,w)-\Phi[f](w+h)+\Phi[f](w)
\ee

Expanding in power series in $w$ equation (\ref{homologeq}) gives

\be
(1-x^2)h_n'(x)+(n-1)\left[(b-a)-(b+a)x\right]h_n(x)=f_n-\Phi[f]_n\  ,\ \ \ n\geq 2
\ee
with the solution 
\be\label{solhomeq}
h_n=J_n[f_n]\ \ \ {\mbox{and\ }}\ \ \Phi[f]_n=\Phi_{n}[f_{n}]
\ee
(see (\ref{opJ}), (\ref{defPhi})).

\subsubsection{Estimates on the solution of the homological equation.}

Let $\rho\in I$, $R>0$, $\delta\in (0,\frac{1}{2})$. Denote by $\Pi$ the Banach space of functions $f(x,w)$ analytic on the polydisk  
$$\Delta_{\rho,R}=\{ (x,w)\, ;\, |x-c|<\rho,\ |w|<R\}$$
and continuous up to the boundary, with the sup-norm, denoted $\|\ \|$. 

Let $\Delta_\delta=\Delta_{\rho{\rm{e}}^{-\delta},R{\rm{e}}^{-\delta}}$ and let $\Pi_\delta$ be the space of functions analytic on the polydisk $\Delta_\delta$, continuous up to the boundary; the sup-norm on this space is denoted by 
$\|\ \|_\delta$.

For $f\in\Pi$ having a zero of order 2 at $w=0$ we have
\be\label{serf}
f(x,w)=\sum_{n\geq 2}f_n(x)w^n=\sum_{n\geq 2,k\geq 0}f_{n,k}w^n(x-c)^k
\ee
with
\be\label{estfnk}
|f_{n,k}|\leq \|f\|R^{-n}\rho^{-k}\ ,\ \ \  \ \ \sup_{D\rho}|f_n(x)|\leq\|f\|R^{-n}
\ee

By (\ref{solhomeq}) we have

\be\label{ser}
\Phi[f](w)=\sum_{n\geq 2}\Phi_n[f_n]w^n\ \ \ \ \ \ \  \ \ \  \ h(x,w)=\sum_{n\geq 2}J_n[f_n](x)w^n
\ee

\begin{Lemma}\label{L3}
If $f\in\Pi$ the series (\ref{ser}) converge absolutely on $\Pi_\delta$ and the sums satisfy the estimates
\be\label{L3i}
\|\Phi[f](w)\|_\delta\leq c_1\delta^{-2}\|f\|
\ee
and
\be\label{L3ii}
\|h\|_\delta\leq c_2\delta^{-2}\|f\|
\ee
with the constants $c_1,c_2$ independent of $\rho\in I$ and $\delta\in (0,\frac{1}{2})$.
\end{Lemma}

{\bf{Proof of Lemma \ref{L3}}}

By Lemma \ref{LemmaPhi} and (\ref{estfnk}) we have
$$\|\Phi[f](w)\|_\delta=\big\| \sum_{n\geq 2}\left(f_{n}(c)+\frac{1}{n}\tilde{\Phi}_n[\tilde{f}_n]\right)w^n\big\|_\delta$$
$$\leq \left[\sum_{n\geq 2}\left(1+\frac{c_5}{n}\delta^{-1}\right){\rm{e}}^{-n\delta}\right]\, \|f\|<c_1 \delta^{-2}\|f\|$$
The estimate (\ref{L3ii}) follows directly from (\ref{Estim1}), (\ref{estfnk}), (\ref{solhomeq}).\qed

\

Denote by $\Pi^0$ the functions $f\in\Pi$ with $f(x,0)=0$  and let
\be\label{norm0}
\|f\|^0=\sup_{\Delta_{\rho,R}}\, \frac{|f(x,w)|}{|w|}
\ee

Note that:
\be\label{Note}
\frac{1}{R}\|f\|\leq \|f\|^0,\ \ |f_{n,k}|\leq \|f\|^0\rho^{-k}R^{-n+1},\ \ \sup_{D_\rho}|f_n(x)|\leq \|f\|^0R^{-n+1}
\ee

\begin{Remark}\label{R4}
For $f\in\Pi^0$ we have
$$\|f\|_\delta\leq \delta^{-2}R\|f\|^0$$
\end{Remark}
This estimate is obtained using (\ref{estfnk}) and (\ref{Note}).

\begin{Remark}\label{R5}
For $f\in\Pi^0$ Lemma \ref{L3} holds in $\|\ \|^0$, namely
$$\|h\|_\delta^0\leq c_2\delta^{-2}\|f\|^0$$
\end{Remark}
The proof follows the same steps as the proof of Lemma \ref{L3}.

\begin{Remark}\label{R6}
Let $0<R_1<R_2$. A standard Cauchy estimate and (\ref{Note}) shows that 
if $\phi(w)$ is a function analytic for $|w|<R_2$, continuous for $|w|\leq R_2$, with a zero of order two at $w=0$ then
$$\max_{|w|\leq R_1}|d_w\phi|\leq \frac{1}{1-R_1/R_2}\, \max_{|w|\leq R_2}\frac{|\phi(w)|}{|w|}$$
\end{Remark}

\begin{Lemma}\label{L4}
There exists a constant $\kappa\geq 2$ such that for any numbers $\beta\geq 4$ and $\delta\in(0,\kappa^{-1})$ the following estimate holds: if $\|f\|^0\leq\delta^\beta$ then
$$\| d_wh\|_\delta\leq c_4\delta\ \ \ {\mbox{and}}\ \ \ \|(1+d_wh)^{-1}\|_\delta\leq c_3$$
\end{Lemma}

{\bf{Proof}}

\z Using Remarks\,\ref{R5} and \ref{R6} we have
$$\|d_wh\|_\delta\leq\, \frac{1}{1-{\rm{e}}^{-\delta/2}}\, \|h\|_{\delta/2}^0\, \leq\, \frac{c_2\delta^{-2}}{1-{\rm{e}}^{-\delta/2}}\, \|f\|^0$$
$$\leq\, c'_2\delta^{-3}\, \|f\|^0\, \leq\, c'_2\delta^{\beta-3}\, \leq c_4\delta$$
Let $\kappa\geq 2$ be large enough so that $c'_3\equiv c_4\kappa^{-1}<1$. Then for $\delta\in(0,\kappa^{-1})$ we have
$$ \|(1+d_wh)^{-1}\|_\delta\leq \left(1-\|d_wh\|\right)^{-1}<(1-c'_3)^{-1}=c_3$$ 
which completes the proof of Lemma\,\ref{L4}.\qed

\begin{Lemma}\label{LD}
There exists a number $\beta_2$ depending only on $c_2$ and $\kappa$ such that for any $\beta\geq\max\{\beta_2,2\}$
and for any $\delta\in(0,\kappa^{-1})$ we have: if $\|f\|^0\leq\delta^\beta$ then
$$\|w+h\|_\delta\leq\, R\, {\rm{e}}^{-\delta/2}$$

In particular, $f\left( x,w+h(x,w)\right)\in\Pi_\delta^0$.

\end{Lemma}

{\bf{Proof}}

\z Using Remark\,\ref{R5} we have
$$\|w+h\|_\delta\leq \, R\, {\rm{e}}^{-\delta}\,\left(1+\|h\|^0_\delta\right)\,\leq\, R\, {\rm{e}}^{-\delta}\,\left(1+c_2\delta^{-2}\|f\|^0\right)\,\leq\, R\, {\rm{e}}^{-\delta}\,\left(1+c_2\delta^{\beta-2}\right)$$
which is less than $R\, {\rm{e}}^{-\delta/2}$ if
$\beta\geq 3+{\ln(2c_2)}/{\ln\kappa}\equiv\beta_2$.\qed

\begin{Corollary}\label{Conseq}
If $(x,w)\in\Delta_\delta$ then
$$\sup_{[w,w+h]}\big| d_wf(x,\cdot )-d_w\Phi[f]\big|\leq\frac{{\rm{e}}^{\delta/4}}{1-{\rm{e}}^{-\delta/4}}\, \max_{|w|\leq R{\rm{e}}^{-\delta/4}}\,|f(x,\cdot)-\Phi[f]\big|$$

\end{Corollary}

\begin{Lemma}\label{LE}

Under the assumptions of Lemma\,\ref{LD} we have
$$\max_{|x-c|\leq\rho{\rm{e}}^{-\delta}}\ \max_{|w|\leq R{\rm{e}}^{-\delta/4}}\  \frac{|f(x,w)-\Phi[f](w)|}{|w|}\,\leq\, c_6\delta^{-2}\, \|f\|^0$$
\end{Lemma}

{\bf{Proof}}

\z A Taylor series expansion in $w$ and Lemma\,\ref{LemmaPhi} give
$$\frac{|f(x,w)-\Phi[f](w)|}{|w|}\,\leq\, \sum_{n\geq 2}|f_n(x)-f_n(c)-\frac{1}{n}T_n[f_n]|\, |w|^{n-1}$$
$$\leq\, \sum_{n\geq 2}\,\left[ |x-c|\, \sup\, |f_n'|\, +\, \frac{1}{n}|T_n[f_n]|\right]\, |w|^{n-1}$$
and Lemma\,\ref{LE} follows from (\ref{L3ii}) and Remark\,\ref{R4}. 
\qed


\begin{Lemma}\label{L5}
Under the assumptions of Lemma\,\ref{LD} the difference (\ref{defD}) satisfies
$$\|\mathcal{R}_0[f]\|^0_\delta\, \leq \, c_7\, \delta^{-5}\left(\|f\|^0\right)^2$$
\end{Lemma}

{\bf{Proof}}

\z We have
$$\mathcal{R}_0[f](x,w)\leq\ |h(x,w)|\ \max_{z\in[w,w+h]}\,|d_z\left(f(x,z)-\Phi[f](z)\right) |$$

$$\|\mathcal{R}_0[f]\|^0_\delta\leq\ \|h\|^0_\delta\, \max_{|x-c|\leq \rho{\rm{e}}^{-\delta}}\ \max_{|z|\leq R{\rm{e}}^{-\delta/2}}\,|d_z\left(f(x,z)-\Phi[f](z)\right) |$$
Remark\,\ref{R6} and Lemma\,{\ref{LE} give
$$\leq\, \frac{\|h\|^0_\delta}{1-{\rm{e}}^{-\delta/4}}\, \max_{|x-c|\leq \rho{\rm{e}}^{-\delta}}\ \max_{|z|\leq R{\rm{e}}^{-\delta/4}}\, \frac{|f(x,z)-\Phi[f](z)|}{|z|}\leq\, \frac{\|h\|^0_\delta}{1-{\rm{e}}^{-\delta/4}}\, c_6\delta^{-2}\|f\|^0$$
Remark\,\ref{R5} completes the proof.\qed

\begin{Lemma}\label{L6}

Let $\kappa$ be as in Lemma\,\ref{L4} and assume $\beta>\max\{\beta_2,5\}$.

Then there exists $\kappa_0>\kappa$ so that if $\delta\in (0,\kappa_0^{-1})$ then the remainder (\ref{defR}) satisfies
$\|\mathcal{R}\|^0_\delta\leq \delta^{-6}\left(\|f\|^0\right)^2$.

\end{Lemma}

The {{proof}} is immediate using Lemmas\,\ref{LD}, \ref{L4}, and \ref{L5}. \qed

\subsubsection{The sequence of parameters.}\label{param}

\

\z 1. Choose $\beta>\max\{\beta_2,5\}$. Fix $\rho_1\in (\rho_{min},\rho_0)\subset I$.

\z 2. Choose $\delta_1>0$ small enough, in the following way.

Define recursively the sequence $\delta_{k+1}=\delta_k^{3/2}$. Therefore $\delta_{k+1}=\delta_1^{(3/2)^{k}}$.

Define $R_{k+1}=R_k{\rm{e}}^{-\delta_k}$ and  $\rho_{k+1}=\rho_k{\rm{e}}^{-\delta_k}$. Therefore
$R_{k+1}=R_1{\rm{e}}^{-\eta_k(\delta_1)}$ where $\eta_k(\delta_1)=\delta_1+\delta_1^{3/2}+\delta_1^{(3/2)^2}+\ldots+\delta_1^{(3/2)^k}$

The sequence of continuous functions $\eta_k$ converges uniformly on the interval $[0,\kappa_0^{-1}]$ to a continuous function $\eta$, which satisfies $\eta(0)=0$ and $\eta(\delta_1)>0$ for all $\delta_1>0$.

Choose $\delta_1$ small enough so that the following three inequalities hold:

\z(i) $\rho_1{\rm{e}}^{-\eta(\delta_1)}>\rho_{min}$.\newline
(ii) ${\rm{e}}^{-\delta_1}+c_2\delta_1^4<1$\newline
(iii) $c_4\delta_1<1$.

\z 3. Choose $\nu\geq\max\{\beta,6\}$. 

{\em{Notation.}} Denote $\Delta_k=\Delta_{\rho_k, R_k}$, and the sup-norm on the polydisk  $\Delta_k$ by $\|\cdot\|_k$.

\z 4. Choose $R_1$ small enough, so that the nonlinear term $f$ of the given equation (\ref{perFuchs}) satisfies $\|f\|^0_1\leq\delta_1^\nu$ (possible since $f(x,\cdot)$ has a zero of order two at $w=0$).

\subsubsection{The iteration.}

We construct recursively a sequence of analytic maps.

Let $f^{[1]}=f$ analytic on $\Delta_1$.
We determine $\phi=\phi^{[1]}$ so that equation (\ref{homologeq}) with $f=f^{[1]}$ has a solution $h^{[1]}$ analytic at both $x=\pm1$. By \S\ref{homeq} the map $H^{[1]}(x,w)=w+h^{[1]}(x,w)$ is an analytic equivalence map between equation 
$\mathcal{E}[f^{[1]}-\phi^{[1]}]$ (see (\ref{corperFuchs})) and equation $\mathcal{E}[f^{[2]}]$ where $f^{[2]}=\mathcal{R}[f^{[1]}]$ 
(see (\ref{defR})). Note that $f^{[2]}_2(x)=0$, in other words, $f^{[2]}(x,\cdot)$ has a zero of order three at $w=0$.

The construction is continued inductively. The second step produces a function $\phi^{[2]}$ and a map $h^{[2]}$ so that 
$\mathcal{E}[f^{[2]}-\phi^{[2]}]$ is equivalent to $\mathcal{E}[f^{[3]}]$ by the analytic map $H^{[2]}(x,w)=w+h^{[2]}(x,w)$ and so on.

\subsubsection{Convergence of the sequence of equivalence maps.}\label{convseqem}

Consider the map $H_k=H^{[1]}\circ H^{[2]}\circ\ldots\circ H^{[k]}$ (where composition is in the $w$-variable only).

\begin{Lemma}\label{R7}

(i) $H^{[k]}(\Delta_{k+1})$ is relatively compact in $\Delta_k$.

(ii) $f^{[k]}$, $\phi^{[k]}$ are holomorphic on $\Delta_{k+1}$ and continuous on $\overline{\Delta_{k+1}}$.

(iii) $H^{[k]}$ is a biholomorphism from $ \Delta_{k+1}$ onto its image.

\end{Lemma}

{\bf{Proof}} \

(i) is proved by induction on $k$. For $k=1$, using (in order)  Lemma\,\ref{L3}, \S\ref{param} points 4 and 3, 
Remark\,\ref{Note}, and \S\ref{param} point 2 we obtain that
$|w+h^{[1]}(x,w)|<R_1$ on $\Delta_2$. The other steps are similar.

(ii) follows by Lemmas\,\ref{L3}, \ref{LD} and \ref {L6}.

(iii) $H^{[k]}$ is one-to-one since assuming $w'+h^{[k]}(x,w')=w''+h^{[k]}(x,w'')$ it follows that
$$|w'-w''|=|h^{[k]}(x,w')-h^{[k]}(x,w'')|\leq |w'-w''|\, \sup_{\Delta_k}|d_wh^{[k]}|$$
which, by Lemma \ref{L4} is less or equal than $|w'-w''|c_4\delta_k$. Since $c_4\delta_k<c_4\delta_0<c_4\delta_1$ by \S\ref{param} point 2, it follows that $|w'-w''|=0$. \qed

\


{\bf{Note}}  that all the domains $\Delta_k$ contain the polydisk $\Delta_e$ of radii $(\rho_1{\rm{e}}^{-\eta(\delta_1)},R_1{\rm{e}}^{-\eta(\delta_1)})$.

\begin{Lemma}\label{serc}
(i) The sequence $\{H_k\}_k$ is Cauchy in the norm $\|\cdot\|^0_e$ on $\Pi^0_e$.

(ii) The limit $H=\lim H_k$ is one-to-one and analytic on $\Delta_e$.

\end{Lemma}

{\bf{Proof}} 

Since $\|f^{[k]}\|^0_k\leq \delta_k^\beta$ (by \S\ref{param} point 3) we may apply Lemma \ref{L4} and obtain that $\|d_wh^{[k]}\|_k\leq c_4\delta_k$ and therefore
\be\label{cstar}
\|d_wH_k\|_k\leq\|1+d_wh^{[1]}\|_2\ldots \|1+d_wh^{[k]}\|_{k+1}\leq (1+c_4\delta_1)\ldots(1+c_4\delta_k)<{\rm{e}}^{c_4\eta(\delta_1)}
\ee

Then $$\|H_k-H_{k+1}\|_e^0=\max\frac{|H^{[k]}(x,w)-H^{[k]}(x,w+h^{[k+1]}(x,w))|}{|w|}$$
$$\leq\|d_wH^{[k]}\|_{k+1}\ \|h^{[k+1]}\|^0_e\leq \, {\rm{e}}^{c_4\eta(\delta_1)}\, c_2\, \delta_{k+1}^{-2}\, \|f^{[k+1]}\|^0_{k+1}\leq\, {\rm{e}}^{c_4\eta(\delta_1)}\, c_2\, \delta_{k+1}^{\beta-2}$$
by Remark \ref{R5} and (\ref{cstar}).

Since the series $\sum_k\delta_{k+1}^{\beta-2}$ converges, using Remark\,\ref{R4}, and noting that all the functions $H_k(x,\cdot)$ are perturbations of the identity the result follows. \qed

\subsubsection{Convergence of the sequence of correction maps.}\label{convcorr}

\begin{Remark}\label{const2n}

Suppose the function $f\in\Pi$ has a zero of order (at least) $p$ at w=0: $f_n(x)=0$ for all $n\leq p-1$ 
(see (\ref{serf})). Then $\phi[f]$ has a zero of order $p$ as well (see  (\ref{ser})), while $h$ and $R[f]$ (see (\ref{defR})) have a zero of order $2p-1$.

As a consequence, $f^{[k]}$ and $\Phi^{[k]}$ have a zero of order $2^{k-1}+1$ for all $k\geq 1$.

\end{Remark}

A sequence of holomorphic functions $\Psi^{[k]}(x,u)$ converging uniformly to the correction $\phi(u)$ of Theorem\,\ref{anf} is now constructed as follows.

Denote by $K^{[n]}$ the $w$-inverse of $H^{[n]}$: $K^{[n]}(x,H^{[n]}(x,w))=w$.

Let $\Psi^{[1]}$ be the first correction of $f$: $\Psi^{[1]}=\phi^{[1]}$. Note that $\phi^{[1]}(u)=\phi_2u^2+O(u^3)$.

The next approximation, $\Psi^{[2]}$, of $\phi$ is obtained as the second correction of $f$: 
the image of $f^{[2]}-\phi^{[2]}$ through the $w$-the inverse of the map $u=H_1(x,w)$ is $f(x,u)-\Psi^{[2]}(x,u)$ where
$$\Psi^{[2]}(x,u)=\Phi^{[1]}(x,u)+\frac{1}{\partial_uK^{[1]}(x,u)}\Phi^{[2]}\left(x,K^{[1]}(x,u)\right)$$ 
$$\ \ \ \ \ \ \ \ \ \ \ \ \ \ \ \  \ \ \ \ \ \ \ \  =\Phi^{[1]}(x,u)+\partial_wH^{[1]}\left(x,K^{[1]}(x,u)\right)\Phi^{[2]}\left(x,K^{[1]}(x,u)\right)$$ 

In general, the inverse of $u=H_n(x,w)$ takes $f^{[n+1]}-\phi^{[n+1]}$ into $f-\Psi^{[n+1]}$ where

$$\Psi^{[n+1]}(x,u)=\Phi^{[1]}(x,u)+\partial_wH^{[1]}\left(x,K_1(x,u)\right)\Phi^{[2]}\left(x,K_1(x,u)\right)$$
$$+\partial_wH^{[2]}\left(x,K_2(x,u)\right)\Phi^{[2]}\left(x,K_2(x,u)\right)+\ldots+\partial_wH^{[n]}\left(x,K_n(x,u)\right)\Phi^{[n+1]}\left(x,K_n(x,u)\right)$$ 
where the following notation was used:
$$K_1=K^{[1]},\ {\mbox{and\ }} K_{j+1}(x,u)=K^{[j+1]}\left(x,K_j(x,u)\right)\ \ {\mbox{for\ }}j\geq 1$$

\begin{Lemma}\label{lem2ton}

Consider the Taylor expansion $\Psi^{[n]}(x,u)=\sum_k\Psi^{[n]}_k(x)u^k$.

(i) We have that $\Psi^{[n]}_k(x)\equiv \Psi^{[n]}_k$ do not depend on $x$ for all $k\leq 2^{n-1}+1$.

(ii) Furthermore, $\Psi^{[n]}_k=\phi_k$ for all $k\leq 2^{n-1}+1$.

(iii) The Taylor coefficients of the expansions in powers of $w$ of map $H_{n+1}$ and of the map 
(\ref{serh}) coincide for all powers of $w$ less or equal to $2^{n-1}+1$.

\end{Lemma}

{\bf{Proof}}

\z (i) Follows by induction using Remark\,\ref{const2n}. 

\z (ii) Follows from the fact that the equation $\mathcal{E}[f-\Psi^{[n]}]$ is analytically equivalent (through the map $u=H_n(x,w)$) to an equation with the nonlinear term $f^{[n]}$ having a zero of order $2^n+1$, and the correction is unique. 

\z (iii) Is shown directly by induction. \qed

As a consequence of Lemmas\,\ref{lem2ton} and \ref{serc} the series (\ref{serh}) converges, and therefore so does (\ref{serphi}).

\subsubsection{The multidimensional case.}\label{ddim}

In the vector-valued case the series are split by degree of homogeneity, then the proof is similar to the one-dimensional case. A slight change is that 
in Lemmas\,\ref{LemmaPhi} and \ref{LemmaJ}, which estimate $\Phi_\bfm[\bfF]$ and $J_\bfm[\bfF]$, 
the saddle points (for large $n$) are given by $c_\bfm$, see (\ref{cm}).

\section{Acknowledgments}
The author would like to express warmest gratitude towards D. Lubinsky for his kind help regarding estimates of Jacobi polynomials.

\section{Appendices}

\subsection{Region containing one regular singularity.}\label{rev1s}

Consider a nonlinear perturbation of a system with a regular singular point 

\be\label{gen1sing}
\frac{d\bfu}{dx}=\frac{1}{x}\, L(x)\bfu+\frac{1}{x}\, \tilde{\bff}(x,\bfu)\ \ \ \ \ \ \bfu\in\CC^d,\ x\in\CC
\ee
where ${L}(x)$ is a matrix analytic at $x=0$ and the nonlinear term $\frac{1}{x}\, \tilde{\bff}(x,\bfu)$ can be allowed to have a first order pole at the singular point $x=0$.

If the eigenvalues of $L(0)$ are nonresonant, after an analytic change of coordinates it can be assumed that $L(x)$ is constant (see \cite{Norm_Form} for more details), hence consider
\be\label{1sing}
\frac{d\bfu}{dx}=\frac{1}{x}\, L\bfu+\frac{1}{x}\, {\bff}(x,\bfu)\ \ \ \ \ \ \bfu\in\CC^d,\ x\in\CC
\ee

It is assumed that $\bff(x,\bfu)$ is analytic for $|\bfu |<\rho$ ($\rho>0$) and  $x$ in a domain $D_r$ which is either a disk centered at the origin: $|x|<r'$, or an annulus $r''<|x|<r'$.\footnote{Allowing the nonlinear part to be singular at $x=0$ accommodates systems corresponding to higher order equations.} Such systems are generically analytically linearizable \cite{Nonln}:

\begin{Theorem}\label{Th1sing}

Assume that the eigenvalues $\mu_1,...,\mu_d$ of the matrix $L$ satisfy the Diophantine
condition (\ref{DioCond}) for all $l\in\ZZ$ if $D_r$ is an annulus (respectively $l\in\NN$ if $D_r$ is a disk), and for all $s\in\{1,...,d\}$, and 
$\mathbf{k}\in\NN^d$ with $ |\mathbf{k}|\geq 2$.

\z Then {{the system (\ref{1sing}) is analytically equivalent to its
linear part}}  $\mathbf{w}'=\frac{1}{x}L\mathbf{w}$ for $x\in{D}_ {\tilde{r}} \subset D_r$ and $|\bfu |<{\tilde{\rho}}<\rho$.

\end{Theorem}

{\bf{Remarks}} 

{\bf{1.}} The domain ${D}_ {\tilde{r}}$ can be made arbitrarily close to $D_r$ if  $\tilde{\rho}$ is small enough.

{\bf{2.}} The analytic equivalence map which is close to identity is unique if no eigenvalue $\mu_j$ is integer.

\subsection{Analyticity of (\ref{hnat-1}) at $x=-1$.}\label{details}
Let {$G$} be the {monodromy matrix} of $Y'=MY$ at $x=-1$: after analytic continuation along a closed loop around $x=-1$ the matrix $Y(x)$ becomes {$AC_{(-1)}Y(x)=Y(x)G$}.

Then the analytic continuation of (\ref{hnat-1}) on a closed loop around $x=-1$ yields

\be
AC_{(-1)}{\bfh_n(x,w)=Y(x)G\int_{-1}^xQ(t)^{-1}G^{-1}Y(t)^{-1}\mathbf{R}_n(t,Y(t)GG^{-1}Y(x)^{-1}\bfw)dt}
\ee
\be\nonumber
=\bfh_n(x,w)
\ee

\z which means that (\ref{hnat-1}) are the coefficients of the unique linearization map of  (\ref{perFuchs}) which is analytic at $x=-1$.


\end{document}